\documentclass[10pt]{article}
\usepackage{latexsym}\usepackage{amsbsy}\usepackage{amssymb}
\usepackage[latin1]{inputenc} 

\newcommand{\ep}{\hspace*{\fill}$\Box$}
\newcommand{\eps}{\varepsilon}
\newcommand{\pr}{{\bf Proof. }}
\newcommand{\ms}{\medskip\\}
\newcommand{\cl}{\mbox{\rm cl}}
\newcommand{\R}{\mathbb R}
\newcommand{\N}{\mathbb N}
\newcommand{\C}{\mathbb C}

\newcommand{\K}{\mathbb K}

\newcommand{\gs}{\ensuremath{{\mathcal G}} }

\newcommand{\es}{\ensuremath{{\mathcal E}} }
\newcommand{\esm}{\ensuremath{{\mathcal E}_M} }

\newcommand{\ns}{\ensuremath{{\mathcal N}} }
\newcommand{\ks}{\ensuremath{{\mathcal K}} }

\newcommand{\comp}{\subset\subset}
\newcommand{\cinfty}{{\cal C}^\infty}
\newtheorem{thr}{\hspace*{-3mm} \bf}[section]
\newcommand{\bt}{\begin{thr} {\bf Theorem. }}
\newcommand{\et}{\end{thr}}
\newcommand{\bp}{\begin{thr} {\bf Proposition. }}
\newcommand{\bc}{\begin{thr} {\bf Corollary. }}
\newcommand{\blem}{\begin{thr} {\bf Lemma. }}
\newcommand{\bex}{\begin{thr} {\bf Example. }\rm} 
\newcommand{\bexs}{\begin{thr} {\bf Examples. }\rm}
\newcommand{\bd}{\begin{thr} {\bf Definition. }}
\newcommand{\beast}{\begin{eqnarray*}}
\newcommand{\eeast}{\end{eqnarray*}}

\newcommand{\al}{\alpha}
\newcommand{\bet}{\beta} 
\newcommand{\ga}{\gamma}
\newcommand{\Om}{\Omega}\newcommand{\Ga}{\Gamma}\newcommand{\om}{\omega}

\newcommand{\vphi}{\varphi}

\newcommand{\D}{{\cal D}}

\newcommand{\pa}{\partial}

\newcommand{\CC}{{\cal C}}

\newcommand{\nn}{\nonumber}
\newcommand{\beq}{ \begin{equation} }\newcommand{\eeq}{\end{equation} }
\newcommand{\bea}{\begin{eqnarray}}\newcommand{\eea}{\end{eqnarray}}
\newcommand{\beas}{\begin{eqnarray*}}\newcommand{\eeas}{\end{eqnarray*}}
\newcommand{\beqs}{\begin{equation*}}\newcommand{\eeqs}{\end{equation*}}

\newcommand{\GL}{\mbox{GL}}\newcommand{\bfs}{\boldsymbol}
\newcommand{\ben}{\begin{enumerate}}\newcommand{\een}{\end{enumerate}}
\newcommand{\ba}{\begin{array}}\newcommand{\ea}{\end{array}}
\newcommand{\brem}{\begin{thr} {\bf Remark. }\rm}
\newcommand{\ethi}{\end{thr}}

\begin{document}

\title{Generalized functions valued in a smooth manifold}
\author{Michael Kunzinger \footnote{Electronic mail: Michael.Kunzinger@univie.ac.at}
  \\ 
  \\{\small Department of Mathematics, University of Vienna}\\
  {\small Strudlhofg. 4, A-1090 Wien, Austria}}
\date{}  
\maketitle
\begin{abstract}
Based on Colombeau's theory of algebras of generalized functions we introduce
the concepts of generalized functions taking values in differentiable manifolds
as well as of generalized vector bundle homomorphisms. 
We study their basic properties, in particular with respect to some new
point value concepts for generalized functions 
and indicate applications of the resulting theory 
in general relativity.
 
\end{abstract}

\section{Introduction}\label{genfmfintro}
An increasing number of applications of Colombeau's theory of algebras of generalized 
functions (\cite{c1,c2,MObook}) 
to questions of a primarily geometric nature, in particular in general relativity
(cf.\ \cite{vickersESI} and the literature cited therein for a survey) have resulted
in a certain restructuring of the theory itself. For its full version,
this task has been pursued in \cite{found, vim}. For the special version
of the theory, which will also provide the framework of the present article, 
the program of ``geometrizing'' the construction has been initiated in
\cite{RD, ndg}. In particular, in \cite{ndg} generalized sections of vector 
bundles were introduced,
a point value description of Colombeau functions on manifolds was given and a number
of consistency results with linear distributional geometry in the sense of 
\cite{marsden}
were derived.

When analyzing singularities in a global context by means of algebras of generalized
functions one inevitably encounters situations where a concept of generalized mappings
{\em valued} in differentiable manifolds is called for. As examples we mention the flow
of a generalized vector field on a manifold or the notion of geodesic of a distributional
spacetime metric. Both notions have been treated in the Colombeau framework, however, due
to the lack of a global notion of generalized mappings between manifolds this treatment
was confined strictly to the local context (cf.\ e.g., 
\cite{DKP, symm, penrose, geo2, vickersESI, genhyp}).

In the present paper we introduce the space $\gs[X,Y]$ of generalized mappings 
defined on the manifold $X$ and taking values in the manifold $Y$, 
adapted to the requirements of applying algebras of generalized functions to global
analysis in the presence of singularities (Section \ref{genmap}). 
The basic idea underlying the construction 
is to consider
nets of smooth maps $u_\eps: X\to Y$ parametrized by a regularization parameter 
$\eps>0$ such that for any compact set $K$ in $X$ and small $\eps$ the values of 
all $u_\eps$ on $K$ remain in a compact subset of $Y$. In addition we will require 
a moderateness condition on the derivatives of the $u_\eps$ (growth on compact
sets bounded by certain inverse powers of $\eps$. This growth restriction has to be
formulated in local charts. An equivalence relation on such nets will be introduced
by employing Riemannian metrics on the manifolds. However, the definition will turn
out to be independent of the chosen metric. Switching between coordinate descriptions
and global formulations in terms of a (but independent of the particular) Riemannian
metric will be a reoccurring pattern throughout this paper.  

Typical singularities that can be modelled by elements
of $\gs[X,Y]$ are jump discontinuities. Singularities of a less tame nature (e.g.,
delta-type singularities) of course have to be handled when differentiating 
elements of $\gs[X,Y]$. In order to accommodate such less pleasant creatures, i.e.,
to develop a suitable notion of tangent map for elements of $\gs[X,Y]$ 
we therefore have to consider generalized
vector bundle homomorphisms (Section \ref{yyyyy})
which will be allowed to behave more singularly in their vector components. 
In this context a new notion of generalized points in vector bundles is introduced
that allows to characterize generalized sections of vector bundles. 
Finally, we show
that the above ``generalized vector bundle points'' are of independent value in that
they allow a direct transfer of the ``pointwise'' nature of smooth tensor fields
to the generalized context.

In the remainder of the introduction we fix some notations from global analysis and
recall those notions and results from \cite{ndg} which will be needed in our
presentation. 

$X$, $Y$ will denote smooth paracompact Hausdorff manifolds 
of dimensions $n$ and $m$, respectively. 
Vector bundles with  base space  $X$  will be denoted by $(E,X,\pi)$ 
(or $(E,X,\pi_X)$)
and for a chart $(V,\vphi)$ in $X$, a vector bundle chart $(V,\Phi)$ 
over $\varphi$ will be written in the form 
($\K=\R$ for real resp.\ $\K=\C$ for complex vector bundles) 
\begin{equation} \label{vbchart}
{\displaystyle
\begin{array}{rcl}
      \Phi:\,\pi^{-1}(V)&\to&\vphi(V)\times\K^{n'}\nn\\
        z&\to&(\vphi(p),\vphi^1(z),\dots,\vphi^{n'}(z))
        \equiv(\vphi(p),\bfs{\vphi}(z))\,, 
\end{array}
}
\end{equation}
where $p=\pi(z)$ and the typical fiber is $\K^{n'}$.
Let $(V_\al,\Phi_\al)_\al$ denote a vector bundle atlas; then we write 
$\Phi_\al\circ\Phi_\beta^{-1}(y,w) = (\vphi_{\al\bet}(y),\bfs{\vphi_{\al\beta}}
(y)w)$, where $\vphi_{\al\bet}:=\vphi_\al\circ\vphi^{-1}_\bet$ is the change 
of chart on the base and $\bfs{\vphi_{\al\beta}}:$
$\vphi_\bet(V_\al\cap V_\beta)\to \GL({n'},\K)$ denotes the transition functions.

For vector bundles $(E,X,\pi_X)$, $(F,Y,\pi_Y)$ we denote by $\mathrm{Hom}(E,F)$
the space of smooth vector bundle homomorphisms from $E$ to $F$. If $f\in 
\mathrm{Hom}(E$, $F)$ we denote by $\underline{f}: X\to Y$ the smooth map
induced by $f$ on the base manifolds. Thus $\pi_Y\circ f = \underline{f}\circ \pi_X$.
For vector bundle charts $(V,\Phi)$ of $E$ and $(W,\Psi)$ of $F$ we 
write the local vector bundle homomorphism 
$f_{\mathrm{\Psi}\mathrm{\Phi}}:=
\mathrm{\Psi}\circ f \circ \mathrm{\Phi}^{-1}: \vphi(V\cap \underline{f}^{-1}(W))
\times \K^{n'} \to \psi(W) \times \K^{m'}$
in the form
\begin{equation}\label{vbhomloc}
f_{\mathrm{\Psi}\mathrm{\Phi}}(x,\xi) = 
(f_{\mathrm{\Psi}\mathrm{\Phi}}^{(1)}(x),f_{\mathrm{\Psi}\mathrm{\Phi}}^{(2)}(x)
\cdot\xi)\,.
\end{equation}

The space of $\cinfty$-sections of $(E,X,\pi)$ is denoted by $\Ga(X,E)$.
As a special case we mention $\Ga(X,T^r_s(X))={\mathcal T}^r_s(X)$, 
the space of $(r,s)$-tensor fields over $X$.

Turning now to notations from Colombeau theory 
we set $I=(0,1]$, 
let ${\cal P}(X)$ be the space of linear differential operators on $X$ and define
\beas
        {\mathcal E}(X) &:=& \cinfty(X)^I\\
        \esm(X)&:=&\{ (u_\eps)_\eps\in{\mathcal E}(X):\ 
        \forall K\subset\subset X,\ \forall P\in{\cal P}(X)\ \exists N\in\N:
  \\&&
        \sup_{p\in K}|Pu_\eps(p)|=O(\eps^{-N})\}
   \\
        \ns(X)&:=& \{ (u_\eps)_\eps\in\esm(X):\ 
        \forall K\subset\subset X,\ \forall q \in\N_0:\
        \sup_{p\in K}|u_\eps(p)|=O(\eps^{q}))\}\,.
\eeas
$\gs(X):= \esm(X)/\ns(X)$ is called the (special) Colombeau algebra on $X$.
$u\in\gs(X)$ if and only if $u\circ\psi_\al\in\gs(\psi_\al(V_\al))$ 
(the local Colombeau algebra on $\psi_\al(V_\al)$)
for all charts $(V_\al,\psi_\al)$. Hence $\gs(X)$ can be identified with the
space of families $(u_\alpha)_\alpha$ with $u_\al \in \gs(\psi_\al(V_\al))$ 
satisfying $u_\al|_{\psi_\al(V_\al\cap V_
\beta)}=u_\beta|_{\psi_\beta( V_\al\cap V_\beta)}\circ\psi_\beta\circ\psi_\al^{-1}$
for all $\al,\beta$ with $V_\al\cap V_\beta\not=\emptyset$. $\cinfty(X)$ is a subalgebra
of the differential algebra $\gs(X)$ and there exist injective sheaf morphisms 
embedding $\D'(X)$, the space of Schwartz distributions on $X$, into $\gs(X)$.

A net $(p_\eps)_\eps \in X^I$ is called compactly supported if there exists some 
$K\comp X$ and some $\eps_0>0$ such that $p_\eps \in K$ for $\eps<\eps_0$. Two
compactly supported points $(p_\eps),\, (q_\eps)_\eps$ are called equivalent,
$(p_\eps)_\eps\sim (q_\eps)_\eps$, 
if $d_h(p_\eps,q_\eps) = O(\eps^m)$ for each $m>0$ for the distance function
$d_h$ induced on $X$ by a Riemannian metric $h$ on $X$. This notion is independent
of the particular $h$. The quotient space $\tilde X_c$ of the set of compactly 
supported points modulo this equivalence relation is called the space of
compactly supported generalized points on $X$. Its elements are written as
$\tilde p = \cl[(p_\eps)_\eps]$. For $u\in \gs(X)$, $u(\tilde p):=
\cl[(u_\eps(p_\eps))_\eps]$ yields a well-defined element of $\ks$, the space
of generalized numbers (see \cite{point}). Any $u\in \gs(X)$ is uniquely determined
by its point values on $\tilde X_c$ (\cite{ndg}, Th.\ 1).

Generalized sections of $E$ are defined as elements of 
$\Gamma_\gs(X,E)$, where (with ${\cal P}(X,E)$ 
the space of linear differential operators on $\Ga(X,E)$))
\beas
        \Gamma_{\esm}(X,E)&:=& \{ (s_\eps)_{\eps\in I}\in \Gamma(X,E)^I :
                \ \forall P\in {\cal P}(X,E)\, \forall K\comp X \, \exists N\in \N:\\
                 &&\sup_{p\in K}\|Pu_\eps(p)\|_h = O(\eps^{-N})\}\\
        \Gamma_\ns(X,E)&:=& \{ (s_\eps)_{\eps\in I}\in \Gamma_{\esm}(X,E) :
                \ \forall K\comp X \, \forall m\in \N:\\
                 &&\sup_{p\in K}\|u_\eps(p)\|_h = O(\eps^{m})\}\,,\ 
                   \mbox{ and, finally}\\
        \Gamma_\gs(X,E)&:=&\Gamma_{\esm}(X,E)/\Gamma_\ns(X,E)\,.
\eeas
Here $\|\,\|_h$ denotes the norm induced on the fibers of $E$ by any Riemannian
metric $h$ on $E$.

The $\cinfty(X)$-module $\Gamma_\gs(X,E)$ is algebraically characterized by
(see \cite{ndg}, Th.\ 4)
\begin{equation}\label{tensorp}
  \Gamma_\gs(X,E)=\gs(X)\otimes_{\cinfty(X)}\Ga(X,E)\,,
\end{equation}

Finally, generalized tensor fields are introduced as generalized sections
of the bundle $T^r_s(X)$. The space of generalized $(r,s)$-tensor fields is
denoted by ${\gs}^r_s(X)$. Algebraically, the following characterizations have
been established in \cite{ndg}, Th.\ 7.
\beas
   \gs^r_s(X)&\cong&L_{\CC^\infty(X)}({\mathfrak X}^*(X)^r,{\mathfrak X}(X)^s;\gs(X))\\
   \gs^r_s(X)&\cong&L_{\gs(X)}(\gs^0_1(X)^r,\gs^1_0(X)^s;\gs(X))\\
   \gs^r_s(X)&\cong&\gs(X)\otimes_{\cinfty(X)} {\cal T}^r_s(X)\,.
\eeas
Here, the first line gives a $\cinfty(X)$-module isomorphism, the second a
$\gs(X)$-isomorphism and the third line is valid in both senses.

\section{Generalized mappings valued in a manifold}\label{genmap}

When generalizing the space of Colombeau generalized functions on the open set
$\Om$ taking ``values'' in $\R^m$ to the manifold setting, the first problem,
namely passing from $\Om$ to $X$ has been dealt with successfully in
\cite{RD, ndg}. Any attempt to replace the range space $\R^m$ by some
manifold $Y$ in a straightforward manner, however, leads into serious troubles due to
the linear structure of $\R^m$ being used in an essential way when
defining moderateness and negligibility, respectively.

The first idea which comes to mind, of course, is to take a local
point of view and to try to formulate moderateness and
negligibility in terms of coordinate representations of
representatives $(u_\eps)_\eps$. Yet, transforming the functions
$u_\eps:\R\to(0,2)$, $u_\eps(x):\equiv\eps$ by the diffeomorphism
$\psi:(0,2)\to(e^{1/2},\infty)$, $\psi(y):= e^{1/y}$ on the range
space shows that
this na\"\i ve approach is bound to fail since it does not even lead to a
well-defined notion of moderateness.

A solution to this problem can be based on generalizing the space
$\gs[\Om,\Om']$ 
introduced in \cite{AB}, 7.3 (denoted there by $\gs_*(\Om,\Om')$) 
to the manifold setting. 
We shortly recall its 
definition. 
Let 
$$
\es[\Om,\Om'] = \{(u_\eps)_\eps \in \cinfty(\Om,\Om')^I
\mid \forall K\comp \Om \exists \eps_0 \exists K'\comp \Om' \forall \eps<\eps_0
\, u_\eps(K) \subseteq K'\}
$$
and denote by $\esm[\Om,\Om']$ the subset of moderate elements of $\es[\Om,\Om']$,
i.e., we set $\esm[\Om,\Om'] = \es[\Om,\Om']\cap \esm(\Om)^m$.
Then $\gs[\Om,\Om'] = \esm[\Om,\Om']/(\ns(\Om)^m)$. 
In \cite{AB}, $\gs[\Om,\Om']$ is called the space 
of Colombeau generalized functions {\em valued} in $\Om'$. Since in this 
terminology the notion ``generalized function valued in $\R$'' would be ambiguous
(referring to both $\gs(\Om)$ and $\gs[\Om,\R]$) we prefer to adopt the term
``c(ompactly)-bounded'' for the generalized functions under consideration.

In defining $\esm[\Om,\Om']$ 
the functions themselves (i.e., derivatives
of order zero) are controlled by a purely topological
condition while, as it will turn out, for the equivalence relation
corresponding to factoring out 
the null ideal, one of the conditions on
the zero-th derivative can be formulated on a manifold 
by means of a distance
function generated by some Riemannian metric.
Thus in both cases, we are led to conditions of intrinsic character
on the manifold. Concerning general derivatives (of order at least 1
in the case of moderateness),
the usual estimates can be transferred to the manifold setting using local
language, due to the chain rule.

There is one technicality to be
observed in the latter process: Fixing a compact subset $K$ on the
manifold $X$ and studying the behavior of $(u_\eps)_\eps$ on $X$
by means of local charts $(V,\varphi)$
and $(W,\psi)$ in $X$ resp.\  $Y$, it may well happen that $K$
contains boundary points of $V$ and/or $u_\eps^{-1}(W)$. Towards the
boundary, however, the functions constituting a chart may tend to
infinity arbitrarily fast. Therefore, we have to restrict our
attention to compact subsets $L$ and $L'$
of $V$ resp.\ $W$ when controlling growth of derivatives in local
terms.

Summarizing, our definitions
eventually will be of partly geometric and partly local
character. The following
proposition suggests the correct way of extracting relevant information
on the moderateness of $(u_\eps)_\eps$ from its local expressions
in charts.

\bp \label{gomomchar}
Let $(u_\eps)_\eps \in \cinfty(\Om,\Om')^I$. The following are
equivalent:
\begin{itemize}
\item[(a)] $(u_\eps)_\eps \in \esm[\Om,\Om']$.
\item[(b)] 
  \begin{itemize}
  \item[(i)] $\forall K\comp\Om
             \ \exists \eps_0>0\  \exists K'\comp \Om' \  \forall 
             \eps<\eps_0:\ u_\eps(K) \subseteq K'$.
  \item[(ii)] 
   $\forall k\in\N$,
   for each chart 
   $(V,\vphi)$
   in $\Om$, each chart $(W,\psi)$ in $\Om'$, each $L\comp V$
   and each $L'\comp W$
   there exists $N\in \N$ with
   $$ \sup\limits_{x\in L\cap u_\eps^{-1}(L')}
      \|D^{(k)}(\psi\circ u_\eps \circ
      \vphi^{-1})(\vphi(x))\|=O(\eps^{-N}).
   $$
  \end{itemize}
\end{itemize}
\et
{\bf Remark.} (i) Here $D^{(k)}$ denotes the operation of
forming
the total derivative of order $k$ and $\|\,\,\|$
is any norm on the respective space of multilinear maps 
(the condition is clearly independent of the choice of $\|\,\,\|$).

(ii) Explicitly, the estimate in (b)(ii) means:
\beas
\lefteqn{
\exists\eps_1>0\ \exists C>0\ \forall\eps<\eps_1\ \forall
x\in
L\cap u_\eps^{-1}(L'):}
\\&&\hphantom{mmmmmmmmm}
\|D^{(k)}(\psi\circ u_\eps \circ
      \vphi^{-1})(\vphi(x))\|
\leq C\eps^{-N}.
\eeas
Therefore, let us agree to treat a supremum over the empty set as
having the value zero in the present context.

\vskip3pt
\noindent\pr (a)$\Rightarrow$(b): (i) holds by definition. Concerning (ii), 
let $k\in \N$ and  choose $N\in \N$, $\eps_1>0$ 
such that $\|D^{(j)}u_\eps(x)\|\le \eps^{-N}$ for $1\le j\le k$, $\eps<\eps_1$ and
all $x\in L$. Then by the chain rule, for $x\in L$ with $u_\eps(x)\in L'$,
$\|D^{(k)}(\psi\circ u_\eps\circ \vphi^{-1})(\varphi(x))\|$ can be 
estimated by $C\eps^{-N}$, where $C$ depends only on 
$\sup_{y\in L'}\|D^{(j_1)}\psi(y)\|$ and 
$\sup_{z\in \vphi(L)}\|D^{(j_2)}(\vphi^{-1})(z)\|$
($1\le j_1,\, j_2\le k$), so the claim follows.

\noindent (b)$\Rightarrow$(a): For estimating derivatives
of order at least 1,
simply set $V=\Om$, $\vphi=\mathrm{id}_\Om$, $L=K$,
$W=\Om'$, $\psi=\mathrm{id}_{\Om'}$ and $L'=K'$
where $K,K'$ are as in the definition of $\es[\Om,\Om']$.
\ep

\vskip6pt
Based on this local characterization we now proceed to introduce 
the concept of moderate mappings between manifolds $X$ and $Y$.
\bd \label{modmapmf} 
The space $\esm[X,Y]$  
of compactly bounded 
(c-bounded)
moderate maps from $X$ to $Y$ is defined as the set
of all $(u_\eps)_\eps \in \cinfty(X,Y)^I$ such that
\begin{itemize}
  \item[(i)] $\forall K\comp\Om\ \exists \eps_0>0\  \exists K'\comp Y \  \forall 
             \eps<\eps_0:\ u_\eps(K) \subseteq K'$.
  \item[(ii)] 
   $\forall k\in\N$,
   for each chart 
   $(V,\vphi)$
   in $X$, each chart $(W,\psi)$ in $Y$, each $L\comp V$
   and each $L'\comp W$
   there exists $N\in \N$  with
   $$\sup\limits_{x\in L\cap u_\eps^{-1}(L')} \|D^{(k)}
   (\psi\circ u_\eps \circ \vphi^{-1})(\vphi(p))\| =O(\eps^{-N}).
   $$
\end{itemize}
\et
\brem \label{modoneatlas}
(i) Membership of $(u_\eps)_\eps$ in $\esm[X,Y]$ 
can equivalently be determined by requiring (ii) merely for charts 
from given atlases ${\cal A}_X = \{(V_\al,\vphi_\al)\mid \al\in A\}$ of $X$ and  
${\cal A}_Y = \{(W_\beta,\psi_\beta)\mid \beta\in B\}$ of $Y$. Indeed, suppose that
(ii) holds for all charts from
${\cal A}_X$, ${\cal A}_Y$ and let $k\in\N$, $L\comp V$, $L'\comp W$, 
with $(V,\vphi)$, $(W,\psi)$ charts in $X$ and $Y$, respectively. 
Choose 
$V_{\al_1},\dots,V_{\al_r}$ covering $L$ and $W_{\beta_1},\dots,W_{\beta_s}$
covering $L'$. We may write $L = \bigcup_{i=1}^r L_i$, $L'=\bigcup_{j=1}^s L_j'$
with $L_i \comp V_{\al_i}$, $L_j'\comp W_{\beta_j}$.
Then for each pair $i,j$ we obtain constants
$N_{ij}$, $\eps_{ij}$ and $C_{ij}$ such that the corresponding
estimate is satisfied on $L_{\al_i}\cap u_\eps^{-1}L_{\bet_j}'$, for
all $\eps<\eps_{ij}$. Setting $N:=\max N_{ij}$
and $\eps_1:=\min \eps_{ij}$,
let $\eps<\eps_1$ and $p\in L$ such that $u_\eps(p)\in L'$.
Suppose that 
$p\in L_{i}$ and $u_\eps(p)\in L'_j$. Then in some neighborhood
of $p$, $\psi\circ u_\eps\circ\vphi^{-1} = \psi\circ \psi_{\beta_{j}}^{-1}
\circ(\psi_{\beta_{j}}\circ u_\eps\circ \vphi_{\al_{i}}^{-1})
\circ \vphi_{\al_{i}}\circ\vphi^{-1}$ and the desired estimate follows from the
fact that any derivative of $\psi\circ \psi_{\beta_{j}}^{-1}$ is uniformly bounded 
on $\psi_{\beta_{j}}(L'_{j})$ and similar for $\vphi_{\al_{i}}\circ\vphi^{-1}$ on $\vphi(L_{i})$. 

\vskip6pt
(ii) Expanding once more the Landau symbol in (ii), we obtain (among
others) $\exists N\in \N\dots\exists\eps_1>0\ \exists C>0
\ \forall\eps<\eps_1\dots{}\le
C\eps^{-N}$. Let us emphasize that, in a certain sense, $N$ and
$\eps_1$ have an intrinsic meaning whereas $C$ depends on the
charts under consideration. The second part of this assertion is
immediate by inspecting the proof of Proposition \ref{gomomchar}; the first part
can be made precise as follows: Given $K\comp\Om$, $l$ and $\eps_1$
can be chosen independently of $V,\varphi,W,\psi,L,L'$ as long as
$L$ is allowed to vary only in $K$. This can be verified by an
argument similar to that in the first part of this remark:
Choosing $K'$ according to (i) and covering $K$ and $K'$ by
the domains $V_i$ resp.\ $W_j$ of (finitely many) suitable charts,
$K$ and $K'$ can be decomposed according to this covering.
Each pair $i,j$ giving rise to constants $N_{ij}$, $\eps_{ij}$,
$C_{ij}$, it is not difficult to show that again $N:=\max N_{ij}$
and $\eps_1:=\min \eps_{ij}$ satisfy (the expanded form of) (ii) of
Definition \ref{modmapmf} for all $V,\varphi,W,\psi,L,L'$as above, provided
$L\subseteq K$, in addition.
\et

From Remark \ref{modoneatlas}\, (i) it follows 
that $(u_\eps)_\eps \in {\cal C}^\infty(X,Y)^I$
is in $\esm[X,Y]$ if and only if $(u_\eps|_V)_\eps \in \esm(V,Y)$ for each $V\in
{\cal V}$ where ${\cal V}$ is any open cover of $X$.

\bd\label{equrel} Two elements $(u_\eps)_\eps$, $(v_\eps)_\eps$ of $\esm[X,Y]$ are 
called equivalent, $(u_\eps)_\eps \sim (v_\eps)_\eps$, 
if the following conditions are satisfied: 
\begin{itemize}
\item[(i)]  For all $K\comp X$, $\sup_{p\in K}d_h(u_\eps(p),v_\eps(p)) \to 0$ 
($\eps\to 0$)
for some (hence every) Riemannian metric $h$ on $Y$.
\item[(ii)] $\forall k\in \N_0\ \forall m\in \N$,
for each chart 
   $(V,\vphi)$
   in $X$, each chart $(W,\psi)$ in $Y$, each $L\comp V$
   and each $L'\comp W$:
$$
\sup\limits_{x\in L\cap u_\eps^{-1}(L')\cap v_\eps^{-1}(L')}\!\!\!\!\!\!\!\!\!\!\!\!\!\!\!\!\!\!\!
\|D^{(k)}(\psi\circ u_\eps\circ \vphi^{-1}
- \psi\circ v_\eps\circ \vphi^{-1})(\vphi(p))\|
=O(\eps^m). 
$$
\end{itemize}
\et
That Definition \ref{equrel}\, (i) is indeed  independent of the 
particular Riemannian metric 
$h$ follows from Definition \ref{modmapmf}\, (i) and Lemma 2 from \cite{ndg}.
To prove 
that also in Definition \ref{equrel}\, (ii) it suffices to consider 
charts from given atlases
(and in several instances to follow)
we shall make use of the following auxiliary result. 
\blem \label{kriegllem}
Let $f:\Om\to \Om'$
be continuously differentiable and suppose that $K\comp \Om$. Then there exists
$C>0$ such that $\|f(x) - f(y)\| \le C\|x-y\|$ for all $x,y\in K$. 

$C$ can be chosen as $C_1\cdot\sup_{z\in L}(\|f(z)\|+\|Df(z)\|)$
where $L$ is any fixed compact neighborhood of $K$ in $\Om$ and
$C_1$ only depends on $L$.
\et
\pr
Choose
a bump function $\chi$ compactly supported in $\Om$ and identical to $1$ in a 
neighborhood of $K$. Then $\tilde f:= \chi f$ can be extended to a smooth function
on all of $\R^n$ and we may take for $C$ the supremum of $\|D\tilde f\|$ on the convex
hull of $K$. The form of $C$ given above follows from the Leibniz
rule.
\ep

\vskip6pt
\brem \label{equivoneatlas}
(i) Equivalence of $(u_\eps)_\eps$, $(v_\eps)_\eps \in \esm[X,Y]$ 
can equivalently be determined by requiring (ii) of Definition \ref{equrel} 
merely for charts 
from given atlases ${\cal A}_X = \{(V_\al,\vphi_\al)\mid \al\in A\}$ of $X$ and  
${\cal A}_Y = \{(W_\beta,\psi_\beta)\mid \beta\in B\}$ of $Y$.
Indeed, suppose that (ii) holds for all charts from
${\cal A}_X$, ${\cal A}_Y$ and let $L\comp V$, $L'\comp W$, 
with $(V,\vphi),\,(W,\psi)$ charts and $k\in\N_0$ and $m\in\N$. 
Choose $V_{\al_1},\dots,V_{\al_r}$ covering $K$ and $W_{\beta_1},\dots,W_{\beta_t}$
covering $L'$. 
We may write $L= \bigcup_{i=1}^r L_i$ and
$L'=\bigcup_{j=1}^t L_j'$ with $L_i \comp V_{\al_i}$, $L_j'\comp W_{\beta_j}$.
For $1\le j\le t$
let $W_j'$ be a neighborhood of $L_j'$ whose closure is compact and 
contained in  $W_{\beta_j}$.
By Definition \ref{equrel}\, (i) there exists
$\eps_1'$ such that
$\sup_{p\in K}d_h(u_\eps(p),v_\eps(p)) < \min_{1\le j \le t}d_h(L_j',\pa W_j')$
for all $\eps<\eps_1'$. 
Now choose $\eps_1''$ as in (the expanded version of the Landau
symbol of)
Definition \ref{equrel},  (ii)
suitable for $k$, $m$, and all sets of data
$(V_{\al_i},\vphi_{\al_i}),(W_{\bet_j},\psi_{\bet_j}),L_i,
\overline{W'_j}$ simultaneously.
Set $\eps_1:=\min(\eps_1',\eps_1'')$.
Then if $p\in L$ and
$\eps<\eps_1$ are such that $u_\eps(p),\, v_\eps(p) \in L'$
there exist $i$ such that $p\in L_i$
and $j$ with $u_\eps(p)$, $v_\eps(p) \in 
W_j'$.

Then 
in a neighborhood of $x:= \vphi(p)$ we have
\beas
\lefteqn{\psi\circ u_\eps\circ\vphi^{-1}-\psi\circ
                                  v_\eps\circ\vphi^{-1}}
\\&&\hphantom{mmm}                                
  =(\psi\circ\psi_{\beta_j})\circ(\psi_{\beta_j}\circ u_\eps
  \circ \vphi_{\al_i}) \circ(\vphi_{\al_i}\circ\vphi^{-1}) 
\\&&\hphantom{mm=(\psi}
  -(\psi\circ\psi_{\beta_j})\circ(\psi_{\beta_j}\circ v_\eps
\circ \vphi_{\al_i})\circ(\vphi_{\al_i}\circ\vphi^{-1}).
\eeas
Hence the claim follows from the chain rule by taking into account Lemma
\ref{kriegllem}.

\vskip6pt
(ii) Again $\eps_1$ occurring in the expanded version of the Landau
symbol in Definition \ref{equrel}\, (ii) has an intrinsic character whereas the
constant $C$ in the respective estimate depends on the
charts under consideration. The proof is similar to that of part
(i) of this remark, compare (ii) of Remark \ref{modoneatlas}.
\et

Similar considerations as in Remark \ref{equivoneatlas} show that 
the relation $\sim$ is transitive. Hence $\sim$ is indeed an equivalence relation. 
Having introduced the basic building blocks of our construction we may now proceed
to defining the space $\gs[X,Y]$ itself.

\bd \label{gmndef} 
The quotient $\gs[X,Y]:=\esm[X,Y]\big/\sim$ is called the space
of compactly bounded (c-bounded) 
Colombeau generalized functions from $X$ into $Y$. 
\et

Due to the c-boundedness of representatives, typical singularities
that can be modelled by elements of $\gs[X,Y]$ are jump
discontinuities along hypersurfaces. When differentiating such
generalized functions (by means of component-wise tangent maps to be
introduced below), however, also delta-type singularities will
enter the picture.

\bexs 
\begin{itemize}
\item[(i)] By Proposition \ref{gomomchar}
and a similar result for the respective null ideal
of $\gs[\Om,\Om']$, for $X=\Om\subseteq \R^n$, $Y=\Om'\subseteq \R^m$ 
we obtain precisely $\gs[\Om,\Om']$ as introduced 
at the beginning of this section.
\item[(ii)] Let $H=\cl[(H_\eps)_\eps] \in \gs(\R)$ such that $H_\eps$ converges
weakly to the Heaviside function. Then $(\exp(i\pi H_\eps))_\eps \in
\esm[\R, S^1]$ may be viewed as a representative of a jump function valued
in $S^1$ (cf.\ also Example \ref{artiex}  (ii) below).
\item[(iii)] Geodesics in generalized spacetimes. Solving geodesic equations
in singular spacetimes modelled locally in Colombeau algebras (cf.\ \cite{geo2}, 
\cite{vickersESI}, \cite{genhyp} and the literature cited therein) yields
representatives $(u_\eps)_\eps$ where each $u_\eps: J\subseteq \R \to X$.
It can be shown that such ``generalized geodesics'' are elements of $\gs[J,X]$.
For details we refer to \cite{curvature}.
\item[(iv)] In \cite{penrose2}, R.\ Penrose introduced a
``discontinuous coordinate transformation'' $T$
transforming the distributional form 
$$
ds^2 = f(x)\delta(u) du^2 - du\,dv + \sum_{i=1}^2 (dx^i)^2
$$
of an impulsive gravitational wave into a continuous form. 
   This transformation was analyzed in the context
of special Colombeau algebras in \cite{penrose}. 
$T$ was recognized as the distributional limit of a
family $(T_\eps)_\eps$ of diffeomorphisms. It follows from (30) in \cite{penrose}
that $T = \cl[(T_\eps)_\eps]$ provides an example of an element of $\gs[X,X]$.
\end{itemize}
\et
\bd \label{equiv0def}
We call two elements $(u_\eps)_\eps$, $(v_\eps)_\eps$
of $\esm[X,Y]$ equivalent of order 
$0$, 
$(u_\eps)_\eps \sim_0 (v_\eps)_\eps$ if
they satisfy Definition \ref{equrel} \, (i) and
(ii) for $k=0$, or, equivalently, (i) and
\begin{itemize}
\item[(ii')] $\forall K\comp\Om\ \forall m\in \N\ \exists \eps_1>0$
$\forall (W,\psi)$ chart in $Y$
$\forall L'\comp W\ \exists C>0$ such that 
$$
|(\psi\circ u_\eps(p)
- \psi\circ v_\eps(p)|\le C \eps^m
$$
for all $\eps<\eps_1$ and all $p\in K$ with $u_\eps(p),\, v_\eps(p) \in L'$.
\end{itemize}
\et
The equivalence of (ii') above and (ii) of Definition \ref{equrel} for $k=0$
(provided $(u_\eps)_\eps,(v_\eps)_\eps\in\esm[X,Y]$)
can be established by the techniques employed in the proofs of Remarks 
\ref{modoneatlas} and \ref{equivoneatlas}.
As in Remark \ref{equivoneatlas} it follows that it suffices to suppose (ii') for
charts from any given atlas ${\cal A}_Y$. The following result gives a 
global characterization of $\sim_0$.
\bt \label{loccharmf} 
Let $(u_\eps)$, $(v_\eps)_\eps \in \esm[X,Y]$. The following are equivalent.
\begin{itemize}
\item[(a)] $(u_\eps)_\eps \sim_0 (v_\eps)_\eps$.
\item[(b)] For some (hence every)
Riemannian metric $h$ on $Y$,
every $m\in \N$ and every $K\comp X$ we have 
$$
\sup_{p\in K}d_h(u_\eps(p),v_\eps(p)) = O(\eps^m) \qquad (\eps\to 0)\,.
$$ 
\end{itemize}
\et
\pr (b)$\Rightarrow$(a):
(i) from Definition \ref{equiv0def} (adopted from Definition \ref{equrel}) is obvious. 
Concerning (ii'), let us first suppose that $L'$ is contained in a 
geodesically convex set $W'$
such that $\overline{W'}\comp W$.
Given $m\in \N$, choose $\eps'>0$, $C'>0$ such that $\sup_{p\in K}d_h(u_\eps(p),v_\eps(p)) 
\le C'\eps^m$ for $\eps<\eps'$. Then for $\eps<\eps'$, $p\in K$ and $u_\eps(p)$, $v_\eps(p)
\in L'$ we have
\begin{equation}\label{dhequ}
d_h(u_\eps(p),v_\eps(p)) =
\int_{a_\eps}^{b_\eps} \|\gamma'_\eps(s)\|_h\,ds
\end{equation}
where $\gamma_\eps:[a_\eps,b_\eps] \to W'$ is the unique geodesic in $W'$
joining $u_\eps(p)$ and $v_\eps(p)$. 
Since $W'$ is relatively compact in $W$ there exists $C''>0$ such that 
$|\xi|\le C''\|T_{\psi(p)}\psi^{-1}(\xi)\|_h$ 
for all $p\in W'$ and all $\xi\in\R^n$ (with $|\,\,|$, say, the Euclidean norm).
Thus 
\beas
|\psi(u_\eps(p)) - 
\psi(v_\eps(p))| &\le& 
\int_{a_\eps}^{b_\eps} |(\psi\circ\gamma_\eps)'(s)| \,ds \\
&\le& C'' \int_{a_\eps}^{b_\eps} \|\gamma'_\eps(s)\|_h \,ds\,. 
\eeas
From this and (\ref{dhequ}) the claim follows. 
To treat the general case where
$L'$ is any compact subset of $W$,
let $W_q$ be a geodesically convex open 
neighborhood of $q$ 
whose closure is compact and contained in $W$,
for each $q\in L'$. Then there exist $W_{q_1},\dots,W_{q_k}$ 
covering $L'$ and we may 
write $L' = \bigcup_{i=1}^k L_i'$ with $L_i'\comp W_{q_i}$. 
Choose $\eps''>0$ such that
$\sup_{p\in K}d_h(u_\eps(p)$, $v_\eps(p))$ $<$ 
$\min_{1\le i \le k}d_h(L_i',\pa W_{q_i})$
for all $\eps<\eps''$. Then if $p\in K$ and
$\eps<\eps''$ are such that $u_\eps(p),\, v_\eps(p) \in L'$
there exists $i$ with $u_\eps(p)$, $v_\eps(p) \in 
W_{q_i}$.
By the first part of the proof 
there exist $\eps'>0$ and $C_i>0$ ($1\le i\le k$)
such that $|\psi\circ u_\eps(p) - \psi
\circ v_\eps(p)| \le C_i \eps^m$ whenever $\eps<\eps'$,
$p\in K$ and 
$u_\eps(p),\, v_\eps(p) \in W_{q_i}$. We may therefore set $\eps_1=\min(\eps',\eps'')$
and $C=\max_{1\le i\le k} C_i$ to establish the claim.

\noindent(a)$\Rightarrow$(b): 
Suppose that (b) is violated. 
Then for some $K\comp X$ we have
\begin{equation}\label{negation}
\exists m_0\in \N \ \forall k\in \N \ \exists p_k\in K \
\exists \eps_k < \frac{1}{k}:\ d_h(u_{\eps_k}(p_k),v_{\eps_k}(p_k))
\ge \eps_k^{m_0}\,.
\end{equation}
Using Definition \ref{modmapmf}\, (i) 
and choosing suitable subsequences if necessary we may suppose
that $p_k \to p \in K$, $u_{\eps_k}(p_k) \to q \in Y$ and $v_{\eps_k}(p_k) \to q$.
Choose a chart $(W,\psi)$ such that $q\in W$. 
For $k$ sufficiently large we
have $u_{\eps_k}(p_k), \,v_{\eps_k}(p_k) \in L' \comp W$
with $L'= \psi^{-1}(\overline{B_r(\psi(q))})$ ($r>0$ suitable).
Hence also the line connecting $\psi(u_{\eps_k}(p_k))$
and $\psi(v_{\eps_k}(p_k))$ is contained in $\psi(L')$. 
We now multiply the Euclidean metric on $\psi(W)$ by
some cut-off function compactly supported in 
$\psi(W)$ and identical
to $1$ in a neighborhood of $\overline{B_r(\psi(q))}$ 
and extend the pullback of the resulting
product under $\psi$ to a Riemannian metric $g$ on $Y$. Then 
from Definition \ref{equiv0def}\, (ii') we conclude
$$
d_g(u_{\eps_k}(p_k),v_{\eps_k}(p_k)) \le |\psi\circ u_{\eps}(p_k)-
\psi\circ v_{\eps_k}(p_k)| = O(\eps^r)
$$
as $\eps\to 0$ for large enough $k$ and arbitrary $r\in \N$. But by Lemma 
2 of \cite{ndg} we have $d_h(u_{\eps_k}(p_k),v_{\eps_k}(p_k)) \le C'
d_g(u_{\eps_k}(p_k),v_{\eps_k}(p_k))$ for some $C'>0$ and all $k$, so we
arrive at a contradiction to (\ref{negation}).
\ep

\vskip6pt
\brem \label{locthcomp}
It is a rather recent development in the theory
of algebras of generalized functions on open subsets $\Om$ of $\R^n$ that a
characterization of the Colombeau-ideal $\ns(\Om)$ as a subspace 
of $\esm(\Om)$ can be given where the asymptotic vanishing of representatives
has only to be supposed in the $0$-th derivative (\cite{found}, Th.\ 13.1
and the remark following it). 
In our present setting this implies 
that for
elements $(u_\eps)_\eps$, $(v_\eps)_\eps$ of $\esm[\Om,\Om']$ in fact 
$(u_\eps)_\eps \sim (v_\eps)_\eps$ is equivalent with 
$(u_\eps)_\eps \sim_0 (v_\eps)_\eps$. It would certainly be desirable to obtain
such an equivalence also in the general case, i.e., for $\gs[X,Y]$. However,
the Taylor argument used in the local case (based in turn on Landau's paper
\cite{Lan}) is not applicable in this situation. In fact, setting 
$f_\eps := \psi\circ u_\eps \circ \vphi^{-1} - \psi\circ v_\eps \circ \vphi^{-1}$
the domain $D_\eps$ of definition of $f_\eps$ may fail to contain
$x+\eps^{m+N}e_i$ for any $x\in D_\eps$, so
an expression of the form $f_\eps(x+\eps^{m+N}e_i)$ (as required
for the proof of \cite{found}, Th.\ 13.1 
may be undefined. This ``minimal size'' of the
domain of definition of $f_\eps$, however, is a necessary requirement for 
estimating
$D^{(k)}f_\eps$ by means of $D^{(k-1)}f_\eps$ and $D^{(k+1)}f_\eps$
(cf.\ \cite{Lan}, Satz 3). Hence there does not seem to be a straightforward 
adaptation of the $\gs[\Om,\Om']$-argument to the manifold setting and we have
to leave this as an open problem. See, however, Remark \ref{domainremedy} below.

\et
Let us now study the action of elements of $\gs[X,Y]$ on generalized points.
To this end we introduce the following terminology. Let $f:(X,g)\to (Y,h)$ be
a smooth map between Riemannian manifolds. For $p\in X$ we denote by 
$\|T_p f\|_{g,h}$ the norm of the linear map $T_p f: (T_p X, \|\,\,\|_g)
\to (T_{f(p)} X, \|\,\,\|_h)$. Then we have
\blem \label{riemnormlem}
Let $(u_\eps)_\eps \in \esm[X,Y]$, $g$ a Riemannian metric on $X$ and 
$h$ a Riemannian metric on $Y$. Then for any $K\comp X$ there exist $\eps_0>0$,
$N\in \N$ and $C>0$ such that
$$
\sup_{p\in K} \|T_p u_\eps\|_{g,h} \le C\eps^{-N}
$$
for each $\eps<\eps_0$.
\et
\pr Immediate from the local expression of $\|T_p u_\eps\|_{g,h}$ and the
definition of moderateness.
\ep

\vskip6pt
\bp \label{gmnpointwelldef} \label{generalized function!composition}
Let $u=\cl[(u_\eps)_\eps] \in \gs[X,Y]$ and $\tilde p \in \widetilde X_c$. Then
$u(\tilde p):=\cl[(u_\eps(p_\eps))]$ is a well-defined element of $\widetilde Y_c$.
\et
\pr Let $(u_\eps)_\eps \sim (v_\eps)_\eps$ and suppose that $\tilde p$ is supported
in $K\comp X$. Then in particular $(u_\eps)_\eps \sim_0 (v_\eps)_\eps$, so by
Theorem \ref{loccharmf} for any Riemannian metric $h$ we have 
$d_h(u_\eps(p_\eps),v_\eps(p_\eps)) \le 
\sup_{p\in K}d_h(u_\eps(p),v_\eps(p)) = O(\eps^m)$
for each $m\in \N$, so $u_\eps(p_\eps) \sim v_\eps(p_\eps)$.

Moreover, we have to show that $(p_\eps)_\eps \sim (p'_\eps)_\eps$ entails 
$(u_\eps(p_\eps))_\eps \sim (u_\eps(p_\eps'))_\eps$. 
We first note that for $\eps$ sufficiently small, $p_\eps$, $p_\eps'$ 
remain in some $K\comp X$. $K$ can be written in the form
$\bigcup_{i=1}^sK_i$ with $K_i\comp C_i$, $C_i$ a connected component of
$X$ ($1\le i \le s$). By the theorem of Nomizu-Ozeki (see, e.g., \cite{hicks})
we may equip $X$ with a geodesically 
complete Riemannian metric $g$. For the distance $d_g$ induced by
$g$ on $X$ we agree upon setting $d_g(p,p') = \infty$ if $p$ and $p'$
do not lie in the same connected component of $X$.
Then for $\eps < \eps_0$ there exists 
a $g$-geodesic $\gamma_\eps: [a_\eps,b_\eps] \to X$ such that 
$$
d_g(p_\eps,p_\eps') = \int_{a_\eps}^{b_\eps} \| \gamma'_\eps(s)\|_g\,ds
$$
(where we have used that for $\eps$ small, $p_\eps$ and $p_\eps'$ necessarily 
lie in the same $C_i$ with $1\le i\le s$ depending on $\eps$). 
Choosing $r_i\in \R_+$, $p_i\in K_i$ such that $\{p\in X \mid d_g(p,p_i) 
\le r_i\} \supseteq K_i$ it follows
that each $\gamma_\eps([a_\eps,b_\eps])$ is contained in 
$K':= \bigcup_{i=1}^s\{p\in X \mid d_g(p,K_i)\le 2r\}$ with $r=\max_{1\le i\le s} r_i$. 
$K'$ is compact by the Hopf-Rinow theorem and our above assumption on $d_g$.
Thus since 
$(u_\eps)_\eps$ is moderate, by Lemma \ref{riemnormlem} there exists $N\in\N$ and for 
any Riemannian metric $h$ on $Y$ there exists $C>0$ such that
\beas
d_h(u_\eps(p_\eps),u_\eps(p'_\eps)) &\le&
\int_{a_\eps}^{b_\eps} \| (u_\eps\circ\gamma_\eps)'(s)\,\|_h\,ds\\
&\le& \int_{a_\eps}^{b_\eps} \|T_{\ga_\eps(s)} u_\eps\|_{g,h}
\|\gamma_\eps'(s)\|_g\,ds\\
&\le& C \eps^{-N} \int_{a_\eps}^{b_\eps} \|  \gamma'_\eps(s)\|_g\,ds
= C\eps^{-N}d_g(p_\eps,p_\eps')\,.   
\eeas
Since $d_g(p_\eps,p_\eps') = O(\eps^m)$ for each $m>0$,
it follows that $(u_\eps(p_\eps))_\eps \sim (u_\eps(p_\eps'))_\eps$, as claimed.
\ep

\vskip6pt
\bp \label{pvthmfval} Let 
$u=\cl[(u_\eps)_\eps]$, $v=\cl[(v_\eps)_\eps]\in \gs[X,Y]$. 
Then $(u_\eps)_\eps \sim_0 (v_\eps)_\eps$ if and only if $u(\tilde p) = v(\tilde p)$
for all $\tilde p \in \widetilde X_c$.
\et
\pr Necessity has already been established in the proof of 
Proposition \ref{gmnpointwelldef}.
Conversely, suppose that $(u_\eps)_\eps \not\sim_0 (v_\eps)_\eps$. Then there exist
$K\comp X$ $m_0\in \N$, $\eps_k<\frac{1}{k}$  and $p_k\in K$ such that 
$d_h(u_{\eps_k}(p_{\eps_k}),v_{\eps_k}(p_{\eps_k})) \ge \eps_k^{m_0}$ for all $k\in \N$.
Setting $p_\eps := p_k$ for $\eps_{k+1}<\eps\le \eps_k$ we obtain an element
$\tilde p = \cl[(p_\eps)_\eps] \in \widetilde X_c$ such that 
$u(\tilde p) \not= v(\tilde p)$.
\ep

\vskip6pt
The following result demonstrates that composition of generalized functions
can be carried out in our present framework. In its formulation (as well as
for some results to follow) we shall make use of the
following condition. For any $(u_\eps)_\eps \in \esm[X,Y]$ we consider
the property
\begin{equation}\label{artificial}
\forall K\comp X \  \exists \eps_0>0
\  \exists (\psi,W) \ \mbox{chart in}\  Y \ \mbox{such that}\  
\overline{\bigcup_{\eps<\eps_0}u_\eps(K)} \comp W\,. 
\end{equation}
It follows from Theorem \ref{loccharmf} that if any representative $(u_\eps)_\eps$ of 
some $u\in \gs[X,Y]$ satisfies (\ref{artificial}) then so does every other.
In this case we say that $u$ itself satisfies (\ref{artificial}).
\bexs \label{artiex}
\begin{itemize}
\item[(i)] (\ref{artificial}) is trivially satisfied for each $u\in \gs[X,Y]$ 
in case $Y$ is an open subset of some $\R^m$.
\item[(ii)] Let
$$
h_\eps(x) := \frac{1+\eps}{2}[\tanh(x) + \tanh(\frac{x}{\eps})]
$$
and $u_\eps(x):= e^{i\pi h_\eps(x)}$. Then $u = \cl[(u_\eps)_\eps] \in \gs[\R,S^1]$
may be viewed as modelling a 
function displaying a jump at $x=0$ from $-i \in S^1$ to $i\in S^1$. Each
$u_\eps$ is onto $S^1$, so no single chart of $S^1$ can contain 
$\bigcup_{\eps<\eps_0} u_\eps(\R)$ for any $\eps_0>0$. 
Nevertheless $u$ satisfies (\ref{artificial}).
\end{itemize}
\et
\bt \label{compmfwell} 
Let 
$u=\cl[(u_\eps)_\eps] \in \gs[X,Y]$, $v=\cl[(v_\eps)_\eps] 
\in \gs[Y,Z]$. Suppose that $(v_\eps)_\eps$ satisfies (\ref{artificial}).
Then $v\circ u := \cl[(v_\eps\circ u_\eps)_\eps]$ is a well-defined
element of $\gs[X,Z]$. 
\et
\pr 
For each $K\comp X$ there exist $K'\comp Y$ and $\eps'_0>0$ such that
$u_\eps(K) \subseteq K'$ and there exist $K''\comp Z$ and $\eps_0''>0$
with $v_\eps(K')\subseteq K''$ for $\eps<\eps_0''$. Hence $v_\eps(u_\eps(K))\subseteq K''$
for $\eps<\min(\eps_0',\eps_0'')$. To finish the proof that $(v_\eps\circ u_\eps)_\eps$
is moderate let $k\in\N$, 
$(V,\vphi)$ a chart in $X$, $L\comp V$, $(U,\zeta)$ a chart in $Z$ and $L''\comp U$.
Choose $\eps_0$ and $L'$ such that $u_\eps(L)\subseteq L'$ for
$\eps<\eps_0$.
Cover $L'$ by charts $(W_i,\psi_i)$ ($1\le i\le r$) and write
$L'=\bigcup_{i=1}^r L_i'$ with $L_i'\comp W_i$.
Choose $\eps_1(\le\eps_0)$ and $N$ according to Definition \ref{modmapmf}
simultaneously for all $1\le k'\le k$, for $(u_\eps)_\eps$ with respect
to all sets of data $(V,\varphi),(W_i,\psi_i),L,L_i$ and
for $(v_\eps)_\eps$ with respect
to all data $(W_i,\psi_i),(U,\zeta),L',L''$.
Then if $\eps<\eps_1$ and $p\in L$ with $v_\eps(u_\eps(p))
\in L''$ there exists $i$ with $u_\eps(p)\in L_i$. 
Thus we may write $D^{(k)}(\zeta\circ v_\eps\circ u_\eps\circ \vphi^{-1})(\vphi(p))$
as $D^{(k)}((\zeta\circ v_\eps\circ \psi_i^{-1})\circ(\psi_i
\circ u_\eps\circ \vphi^{-1}))(\vphi(p))$ and the claim follows from the chain rule.

To show that $v\circ u$ is well-defined let $(u_\eps)_\eps \sim (u_\eps')_\eps$
Then by Proposition \ref{pvthmfval} for any $\tilde p = 
\cl[(p_\eps)_\eps] \in \widetilde X_c$ we have $(u_\eps(p_\eps))_\eps
\sim (u_\eps'(p_\eps))_\eps$ and, consequently, $(v_\eps(u_\eps(p_\eps)))_\eps
\sim (v_\eps(u_\eps'(p_\eps)))_\eps$. Thus $(v_\eps\circ u_\eps)_\eps \sim_0
(v_\eps\circ u_\eps')_\eps$, again by Proposition \ref{pvthmfval}. In
particular, (i) of Definition \ref{equrel} is satisfied.

To establish (ii) let $k\in\N_0$,
$(V,\vphi)$ a chart in $X$, $(U,\zeta)$ a chart in $Z$,
$L\comp V$, $L''\comp U$. Choose
$L'\comp Y$ such that $u_\eps(L)\cup u_\eps'(L) \subseteq L'$
for $\eps$ small. We can write $L'=\bigcup_{i=1}^r L_i'$, with
$L_i'\comp \overline{W_i'}\comp W_i$, $(W_i,\psi_i)$ charts in $Y$
and $W'_i$ respective open neighborhoods of $L'_i$.
Choose a chart $(U_1,\zeta_1)$ as
provided by (\ref{artificial}) for $(v_\eps)_\eps$
and the compact set $K':=\bigcup_{i=1}^r \overline{W_i'}$, yielding
$K'':=\overline{\,\bigcup_{\eps<\eps_2} v_\eps(K')}\comp U_1$
for some positive $\eps_2$.
Now suppose $p\in L$ and $v_\eps(u_\eps(p)),\, 
v_\eps(u_\eps'(p)) \in L''$.
For each fixed $\eps$ sufficiently small we have 
that $u_\eps(p)$, $u_\eps'(p)$ both are contained in the same $W_i'$
for some $i$. Hence
\beas
\lefteqn {D^{(k)}(\zeta_1\circ v_\eps\circ u_\eps\circ \vphi^{-1}
- \zeta_1\circ v_\eps\circ u_\eps'\circ \vphi^{-1}))(\vphi(p))}\\
&&= D^{(k)}((\zeta_1\circ v_\eps\circ \psi_i^{-1})\circ(\psi_i\circ
u_\eps\circ \vphi^{-1})\\
&&\hphantom{mmmmmmmmmmm}- (\zeta_1\circ v_\eps\circ  \psi_i^{-1})\circ(\psi_i\circ
u_\eps'\circ \vphi^{-1}))(\vphi(p)).
\eeas
This
last expression can be estimated using the chain rule, moderateness
of $(v_\eps)_\eps$, equivalence of $(u_\eps)_\eps$,
$(u_\eps')_\eps$, Lemma \ref{kriegllem} and moderateness of $(u'_\eps)_\eps$.
Note that for each $i$, Lemma \ref{kriegllem} is to be applied to the
compact subset $\psi_i(\overline{W'_i})$ of 
the domain of definition 
$\psi_{i}(W_{i}\cap v_\eps^{-1}(U_1))$ of 
$\zeta_1\circ v_\eps\circ \psi_{i}^{-1}$ 
thus also the constants produced by this part of the argument are in
fact independent of $\eps$. Since all derivatives of 
$\zeta\circ\zeta_1^{-1}$ are bounded on the compact subset
$\zeta_1(L''\cap K'')$ of $\zeta_1(U\cap U_1)$ analogous estimates
follow (involving one more application of Lemma \ref{kriegllem})
with $\zeta$ replacing $\zeta_1$ in the above. Hence
$(v_\eps\circ u_\eps)_\eps \sim (v_\eps\circ u_\eps')_\eps$. 

Finally, if 
$(v_\eps)_\eps \sim (v_\eps')_\eps$ then 
$(v_\eps\circ u_\eps)_\eps \sim_0 (v_\eps'\circ u_\eps)_\eps$ is immediate
from property (i) in Definition \ref{modmapmf} of $(u_\eps)_\eps$.
Let $p\in L$ such that $v_\eps(u_\eps(p))$, $v_\eps'(u_\eps(p)) \in L''$.
As above, for small $\eps$ we have $u_\eps(p)\in L_i'$, so
the estimate for
$D^{(k)}((\zeta\circ v_\eps\circ \psi_i^{-1})\circ(\psi_i\circ
u_\eps\circ \vphi^{-1})
- (\zeta\circ v_\eps'\circ  \psi_i^{-1})\circ(\psi_i\circ
u_\eps\circ \vphi^{-1}))(\vphi(p))$ follows immediately from
the equivalence of $(v_\eps)_\eps$ and $(v_\eps')_\eps$.
\ep

\vskip6pt
The map $\sigma: f \mapsto \cl[(f)_\eps]$ provides a canonical embedding of 
$\cinfty(Y,Z)$ into $\gs[Y,Z]$. Composition with smooth functions can
be carried out unrestrictedly:
\bc \label{compsecgxy}
\begin{itemize}
\item[(i)] Let $u\in \cinfty(X,Y)$, $v=\cl[(v_\eps)_\eps]\in \gs[Y,Z]$. Then
$v\circ u := \cl[(v_\eps\circ u)_\eps]$ is a well-defined element of $\gs[X,Z]$.
\item[(ii)] Let $u=\cl[(u_\eps)_\eps]\in \gs[X,Y]$, $v\in \cinfty(Y,Z)$. Then
$v\circ u := \cl[(v\circ u_\eps)_\eps]$ is a well-defined element of $\gs[X,Z]$.
\end{itemize}
\et
\pr The proof of (i) actually is a slimmed-down version of that of Theorem 
\ref{compmfwell}, with the (hard) part dealing with
$(u_\eps)_\eps \sim (u_\eps')_\eps$ simply omitted.
Concerning (ii), we also may proceed along the lines
of proof of Theorem \ref{compmfwell}, noting that if 
$(u_\eps)_\eps \sim (u_\eps')_\eps$ we may estimate terms of the form
$$
\|D^{(k)}((\zeta\circ f\circ \psi^{-1})\circ (\psi \circ u_\eps \circ \vphi^{-1})
- (\zeta\circ f\circ \psi^{-1})\circ (\psi \circ u_\eps' \circ \vphi^{-1}))(\vphi(p))\|
$$
directly by Lemma \ref{kriegllem} since
both $\psi \circ u_\eps \circ \vphi^{-1}$ and $\psi \circ u_\eps'
\circ \vphi^{-1}$ map $L(\comp X)$ into a fixed compact subset of the
domain of definition of $\zeta\circ f\circ \psi^{-1}$ for $\eps$
small.\ep

\vskip6pt
\brem \label{domainremedy}
Condition (\ref{artificial}) can also be employed to remedy the domain problem
discussed in Remark \ref{locthcomp}. 
Thus if $(u_\eps)_ \eps, (v_\eps)_\eps \in \esm[X,Y]$
and $(u_\eps)_\eps$ satisfies (\ref{artificial}) then $(u_\eps)_\eps \sim_0 (v_\eps)_\eps$
implies $(u_\eps)_\eps \sim (v_\eps)_\eps$. In fact, by Theorem 
\ref{loccharmf}, (\ref{artificial}) 
in this case is satisfied simultaneously for $(u_\eps)_\eps$ and $(v_\eps)_\eps$,
so the Taylor argument indicated in Remark \ref{locthcomp} yields the claim.
\et

\section{Generalized vector bundle homomorphisms}\label{yyyyy}
In 
the present section 
we introduce the notion of generalized vector bundle homomorphisms. Among others
this notion will allow us to treat tangent mappings of elements of $\gs[X,Y]$.
\bd \label{homgdef} 
Let ${\esm}^{\mathrm{VB}}[E,F]$ 
be the set of all nets $(u_\eps)_\eps$ $\in$
$\mathrm{Hom}(E,F)^I$ satisfying 
\begin{itemize}
\item[(i)] $(\underline{u_\eps})_\eps \in \esm[X,Y]$.
\item[(ii)] $\forall k\in \N_0\  
\forall (V,\Phi)$ 
vector bundle chart in $E$, 
$\forall (W,\Psi)$ vector bundle chart in $F$,
$\forall L\comp V\ 
\forall L'\comp W\ \exists N\in \N\ \exists \eps_1>0\ 
\exists C>0$ with 
$$
\|D^{(k)}
(u_{\eps \mathrm{\Psi}\mathrm{\Phi}}^{(2)}(\vphi(p)))\|
\le C\eps^{-N}
$$ 
for all $\eps<\eps_1$ and all $p\in L\cap\underline{u_\eps}^{-1}(L')$.
Here $\|\,.\,\|$ denotes any matrix norm. 
\end{itemize}
\et
Again it suffices to require (ii) merely for charts from arbitrary vector bundle atlases
of $E$ resp.\ $F$. The proof is entirely analogous to that of Remark \ref{modoneatlas},
taking into account the general form (\ref{vbhomloc}) of local vector bundle 
homomorphisms.
\bd \label{homgequ}
$(u_\eps)_\eps$, $(v_\eps)_\eps \in {\esm}^{\mathrm{VB}}[E,F]$
are called $vb$-equivalent,  
denoted by $(u_\eps)_\eps \sim_{vb} (v_\eps)_\eps$, 
if
\begin{itemize}
\item[(i)] $(\underline{u_\eps})_\eps \sim (\underline{v_\eps})_\eps$ in
$\esm[X,Y]$.
\item[(ii)] $\forall k\in \N_0\ \forall m\in \N\ \forall (V,\Phi)$ 
vector bundle chart in
$E$, $\forall (W,\Psi)$ vector bundle chart in $F$,
$\forall L\comp V\ \forall L'\comp W
\ \exists \eps_1>0\ \exists C>0$ such that:
$$
\|D^{(k)}(u_{\eps \mathrm{\Psi}\mathrm{\Phi}}^{(2)}
-v_{\eps \mathrm{\Psi}\mathrm{\Phi}}^{(2)})(\vphi(p))\|
\le C\eps^{m}
$$  
for all $\eps<\eps_1$ and all $p\in L\cap\underline{u_\eps}^{-1}(L')
\cap\underline{v_\eps}^{-1}(L')$. 
\end{itemize}
$(u_\eps)_\eps$, $(v_\eps)_\eps \in {\esm}^{\mathrm{VB}}[E,F]$  are called 
$vb$-$0$-equivalent, denoted by $(u_\eps)_\eps \sim_{vb0} (v_\eps)_\eps$, if 
$(\underline{u_\eps})_\eps \sim_0 (\underline{v_\eps})_\eps$ and (ii) 
above holds for $k=0$. 
We set $\mathrm{Hom}_{\gs}[E,F] := {\esm}^{\mathrm{VB}}[E,F]\big/\sim_{vb}$.
\et 
To see that condition (ii) of Definition \ref{homgequ} is well behaved under 
changing vector bundle charts
(and thereby that (ii) only has to be required for arbitrary vector bundle atlases)
we may again proceed analogously to Remark \ref{equivoneatlas}. Thus we only have to
consider the following local situation:
let $g$, $h$ be local vector bundle homomorphisms, abbreviate 
$u_{\eps \mathrm{\Psi}\mathrm{\Phi}}$ by $u_\eps$ (and similar for $v_\eps$)
and note that $(h\circ u_\eps \circ g)^{(2)} = (x,v) \to 
h^{(2)}(u_\eps^{(1)}(g^{(1)}(x)))u_\eps^{(2)}(g^{(1)}(x))g^{(2)}(x)v$. Hence
\beas
\lefteqn{\|(h\circ u_\eps \circ g)^{(2)}(x) - (h\circ v_\eps \circ g)^{(2)}(x)\|
}\\
&& \le \|h^{(2)}(u_\eps^{(1)}(g^{(1)}(x))) - h^{(2)}(v_\eps^{(1)}(g^{(1)}(x)))\|
\|u_\eps^{(2)}(g^{(1)}(x))g^{(2)}(x)\| \\
&&\hspace{1.7em} + \|h^{(2)}(v_\eps^{(1)}(g^{(1)}(x)))\| \|u_\eps^{(2)}(g^{(1)}(x))
- v_\eps^{(2)}(g^{(1)}(x))\| \|g^{(2)}(x)\|.
\eeas
Here the first summand can be estimated using Lemma \ref{kriegllem}
while an estimate of the second one follows immediately from 
moderateness of $(v_\eps)_\eps$ and Definition \ref{homgequ} for $u_\eps^{(2)}$,
$v_\eps^{(2)}$. Derivatives can be treated analogously.

In the case of local vector bundle homomorphisms $u_\eps, v_\eps: 
\Om\times \R^{n'} \to \Om'\times \R^{m'}$ it follows 
that 
$(u_\eps)_\eps \in {\esm}^{\mathrm{VB}}(\Om\times \R^{n'},\Om'\times \R^{m'})$
if and only if $(u^{(1)}_\eps)_\eps \in \esm[\Om,\Om']$ 
and $(u^{(2)}_\eps)_\eps$ is moderate as a
map from $\Om$ into 
$\R^{{m'}^2}$ and $(u_\eps)_\eps \sim_{vb}
(v_\eps)_\eps$ if and only if $u^{(i)}_\eps - v_\eps^{(i)}$ is negligible in the local
sense ($i=1,2$).

For $u\in \mathrm{Hom}_{\gs}[E,F]$ the definition $\underline{u} := 
\cl[(\underline{u}_\eps)_\eps]$
yields a well-defined element of $\gs[X,Y]$ which is uniquely characterized by
the property $\underline{u}\circ\pi_X = \pi_Y\circ u$ (with $\pi_X$, $\pi_Y$ the
bundle projections).

\bd\label{tangmapdef}
For any $u\in \gs[X,Y]$ we define the tangent map 
$$
Tu := \cl[(Tu_\eps)_\eps] \in \mathrm{Hom}_{\gs}(TX,TY)\,.
$$
\et
That $Tu$ is well-defined, i.e., that $(u_\eps)_\eps \sim (v_\eps)_\eps$ entails
$(Tu_\eps)_\eps \sim_{vb} (Tv_\eps)_\eps$ is immediate from the definitions of $\sim$
and $\sim_{vb}$.
\bd \label{vbpoints} 
On any vector bundle $(E,X,\pi)$ we consider the set of $(e_\eps)_\eps \in E^I$
such that 
\begin{itemize}
\item[(i)] $\exists K\comp X$ $\exists \eps_0>0$ such that 
$(\pi(e_\eps))_\eps \in K$ for all $\eps<\eps_0$.
\item[(ii)] For some Riemannian metric $h$ on $E$ inducing the norm $\|\,\,\|_h$
on the fibers of $E$ we have: 
$\exists \eps_1>0$ $\exists N\in \N$ $\exists C>0$ such that
$$
\|e_\eps\|_h \le C\eps^{-N}
$$
for all $\eps<\eps_1$. 
(A different choice of $h$ in this condition only influences $C$.)
\end{itemize}
On this set of $vb$-moderate generalized points 
we 
introduce an equivalence relation $\sim_{vb}$ 
by calling two elements 
$(e_\eps)_\eps$, $(e'_\eps)_\eps$ equivalent if
\begin{itemize}
\item[(iii)] $(\pi(e_\eps))_\eps \sim (\pi(e'_\eps))_\eps$ in $X^I$
(i.e., $d_g(p_\eps,q_\eps)=O(\eps^s)$ for each $s$ and one (hence every)
Riemannian metric $g$ on $X$).
\item[(iv)] $\forall m\in \N 
\ \forall 
(W,\Psi)$ vector bundle chart in $E$ $\forall L'\comp W
\ \exists \eps_1>0\ \exists C>0$ such that
$$
|\bfs{\psi}e_\eps - \bfs{\psi}e'_\eps| \le C\eps^m
$$
for all $\eps<\eps_1$ 
whenever both $\pi_X(e_\eps)$ and $\pi_X(e'_\eps)$ lie in $L'$.
(Here $\Psi=\big( z_p\mapsto
(\psi(p),\bfs{\psi}(z))\big)$, cf.\ (\ref{vbchart}).)
\end{itemize}
The set of equivalence classes is denoted by $E_c^{\sim_{vb}}$. 
\et
By Lemma \ref{kriegllem} 
it suffices to require (iv) 
for charts from any given vector bundle atlas.
Elements of
$E_c^{\sim_{vb}}$ will be written in the form $\tilde e = \cl[(e_\eps)_\eps]$.
$\pi(\tilde e) := \cl[(\pi(e_\eps))_\eps]\in \tilde X_c$ is called base point of 
$\tilde e$. 
For any $\tilde p\in \tilde X_c$ we set
$$
(E_c^{\sim_{vb}})_{\tilde p} := \{\tilde e \in E_c^{\sim_{vb}} \mid \pi(\tilde e)
= \tilde p\}\,.
$$
Let $\tilde e  = \cl[(e_\eps)_\eps]\in (E_c^{\sim_{vb}})_{\tilde p}$ 
with $\tilde p = \cl[(p_\eps)_\eps]$. 
Then there exists a representative $({e'}_\eps)_\eps$ of $\tilde e$ and 
$\eps_0>0$ such that $\pi(e_\eps') = p_\eps$ for all $\eps<\eps_0$. In fact, 
there exists $K\comp X$ such that $p_\eps \in K$ for $\eps$ small and we may choose
vector bundle charts $(\Psi_i,V_i)$ ($1\le i\le k$) such that the $V_i$ cover $K$.
Write $K=\bigcup_{i=1}^k K_i$ with $K_i \comp V_i$ and let $V_i'$ be a neighborhood
of $K_i$ whose compact closure is contained in $V_i$.
Then there is $\eps_1>0$ such that for every $\eps<\eps_1$ there exists $i$
with $\pi(e_\eps),\, p_\eps \in \pi(V_i')$. 
For $\eps<\eps_1$ 
we set $e_\eps' = \Psi_i^{-1}(\psi_i(p_\eps),
\bfs{\psi}_i(e_\eps))$. Then it follows directly from the 
definition of $\sim_{vb}$ that $\tilde e' = \tilde e$, while $\pi(e_\eps') = p_\eps$
for $\eps$ small holds by construction.

This observation paves the way to endowing $(E_c^{\sim_{vb}})_{\tilde p}$
with the structure of a $\ks$-module: for $\tilde e$, $\tilde e'
\in (E_c^{\sim_{vb}})_{\tilde p}$ and $\tilde r \in \ks$ choose 
representatives $(e_\eps)_\eps,\,(e'_\eps)$ as above and define
$\tilde e + \tilde r \tilde e' := \cl[(e_\eps+e_\eps')_\eps]$. 
Again it follows directly from the definition
of $\sim_{vb}$ that $\tilde e + \tilde r \tilde e'$ is well-defined, i.e.,
it depends exclusively on $\tilde e$, $\tilde e'$ and $\tilde r$.

Our next aim is to show that $E_c^{\sim_{vb}}$ provides the appropriate concept 
of point values of generalized sections and generalized vector bundle homomorphisms.
To this end we first prove the following strengthening of \cite{ndg}, Lemma 3.
\blem\label{pvmfequchar}
Let $\tilde p = \cl[(p_\eps)_\eps]$, $\tilde p' = \cl[(p'_\eps)_\eps] \in
\tilde X_c$. Then $\tilde p = \tilde p'$ if and only if 
\begin{itemize}
\item[(i)] $d_h(p_\eps,p_\eps') \to 0$ for one (each) Riemannian metric $h$ on $X$. 
\item[(ii)] For each $m\in \N$, each chart $(V,\vphi)$ of $X$ and each
$L\comp V$ there exist $\eps_1>0$ and $C>0$ with
$$
|\vphi(p_\eps) - \vphi(p_\eps')| \le C\eps^m
$$
whenever $\eps<\eps_1$ and both $p_\eps$, $p_\eps'$ lie in $L$.
\end{itemize}
\et
\pr This is just an appropriately slimmed-down version of the proof of 
Theorem \ref{loccharmf}.
\ep

\vskip6pt
\bp \label{pvwellvb}
\begin{itemize}
\item[(i)] Let $u\in \Gamma_{\gs}(X,E)$ and $\tilde p \in \tilde X_c$. Then
$u(\tilde p) := \cl[(u_\eps(p_\eps))_\eps]$ 
is a well-defined element of $E_c^{\sim_{vb}}$.
\item[(ii)] Let $v\in
\mathrm{Hom}_{\gs}[E,F]$ such that $(\underline{v_\eps})_\eps$ satisfies
(\ref{artificial}) (for charts on $Y$ induced by vector bundle charts on $F$)
and let $\tilde e \in E_c^{\sim_{vb}}$. Then $v(\tilde e)
:= \cl[(v_\eps(e_\eps))_\eps]$ is a well-defined element of $F_c^{\sim_{vb}}$. 
\end{itemize}
\et
\pr (i) Clearly $(u_\eps(p_\eps))_\eps$ is $vb$-moderate. To show that 
$u(\tilde p)$ is well-defined we first claim that $(u_\eps)_\eps
\in \Gamma_{\ns}(X,E)$ implies $(u_\eps(p_\eps))_\eps \sim_{vb}
(0_{p_\eps})_\eps$. In fact, (iii) of Definition \ref{vbpoints} is 
automatically satisfied
and 
(iv) follows since 
$$
\|u_\eps(p_\eps)\|_h \le \sup_{p\in K}
\|u_\eps(p)\|_h = O(\eps^m)
$$ 
for every $m\in \N$ where $K\comp X$ is such that
$p_\eps \in K$ for $\eps$ small. Supposing now that $(p_\eps) \sim (p_\eps')$
we have to show that $(u_\eps(p_\eps))_\eps \sim_{vb} (u_\eps(p_\eps'))_\eps$.
In this case (iii) of Definition \ref{vbpoints} holds by assumption. 
To prove property (iv), 
cover $K$ by vector bundle charts $(V_i,\Phi_i)$ ($1\le i\le k$) such that 
$K=\bigcup_{i=1}^k K_i$, $K_i\comp V_i'\comp V_i$ and such that
in addition $\vphi_i(V_i')$
is convex for each $i$. There exists $\eps_1>0$ such that for each $\eps<\eps_1$
there is $i$ with $p_\eps$, $p_\eps'\in
V_i'$. For such $\eps$ we have
\beas
\lefteqn{\hspace*{-1em}|\bfs{\vphi}_i(u_\eps(p_\eps)) - \bfs{\vphi}_i(u_\eps(p_\eps'))|
}\\[.2em]
&&= |\bfs{\vphi}_i\circ u_\eps \circ 
\vphi_i^{-1}(\vphi_i(p_\eps)) 
- \bfs{\vphi}_i\circ u_\eps \circ 
\vphi_i^{-1}(\vphi_i(p_\eps'))| \\[.2em]
&&\le\max_{1\le i\le k}\,\sup_{x\in \vphi_i(V_i')}
\|D(\bfs{\vphi}_i\circ u_\eps\circ\vphi_i^{-1})(x)\|
|\vphi_i(p_\eps)-\vphi_i(p_\eps')|.
\eeas
Here the first factor is bounded by some $\eps^{-N}$ since $(u_\eps)_\eps$ is
moderate and the second one is $O(\eps^m)$ for each $m$ by 
Lemma \ref{pvmfequchar}.

(ii) To show $vb$-moderateness of $(v_\eps(e_\eps))$ we first note that
$\pi_F(v_\eps(e_\eps))_\eps$  $=$ $(\underline{v_\eps}(\pi_X(e_\eps)))_\eps$ is 
compactly supported, which gives (i) of Definition \ref{vbpoints}. To see
Definition \ref{vbpoints}\, (ii), choose $K\comp X$ containing $p_\eps:=\pi_X(e_\eps)$ 
for $\eps$ small.
Then $K$ is covered by finitely many vector bundle
charts $(V_i,\Phi_i)$. Writing $K$ as a union of compact sets as above, an
estimate of the form
$$
|v_{\eps \Psi\Phi_i }^{(2)}(\vphi_i(p_\eps)) \bfs{\varphi}_i(e_\eps)| \le 
\|v_{\eps \Psi\Phi_i }^{(2)}(\vphi_i(p_\eps))\| 
|\bfs{\varphi}_i (e_\eps)| \le \eps^{-N}
$$
for some $N$ (using the moderateness assumptions on $(v_\eps)_\eps$ and 
$(e_\eps)_\eps$) yields the claim. 

Suppose that $(e_\eps)_\eps \sim_{vb} (e_\eps')_\eps$ are
supported in $K$ (i.e., $p_\eps \in K$ for $\eps$ small
and analogously for $p_\eps'$). We claim that $(v_\eps(e_\eps))_\eps \sim_{vb}
(v_\eps(e_\eps'))_\eps$. Property (iii) of Definition \ref{vbpoints} is immediate from
Proposition \ref{gmnpointwelldef}, applied to $\underline{v}$. To show 
(iv) we write 
$K=\bigcup_{i=1}^r K_i$, $K_i\comp \overline{V_i'}\comp V_i$, $(V_i,\Phi_i)$ vector
bundle charts in $E$ and $V_i'$ an open neighborhood of $K_i$. 
Let $K_1:= \bigcup_{i=1}^r
\overline{V_i'} \comp \bigcup_{i=1}^r V_i$ and (using (\ref{artificial})) choose a 
vector bundle chart $(\Psi,W)$ in $Y$ 
such that 
$\overline{\bigcup_{\eps<\eps_1} \underline{v_\eps}(K_1)} 
\comp W$ for some $\eps_1>0$. 
By Definition \ref{vbpoints}\, (iii), for each sufficiently small $\eps$ there exists
$i$ with $p_\eps,\, p_\eps' \in V'_i$.
Now we note that $(\Psi\circ v_\eps \circ 
\Phi_i^{-1})(\Phi_i(e_\eps)) = 
(\psi\circ\underline{v_\eps}\circ\vphi_i^{-1}(\vphi_i(p_\eps)), 
v_{\eps \Psi\Phi_i}^{(2)}(\vphi_i(p_\eps))\cdot
\bfs{\varphi}_i(e_\eps))$. Hence
\beas
\lefteqn{\hspace*{-2em}|\bfs{\psi}(v_\eps(e_\eps))- \bfs{\psi}(v_\eps(e_\eps'))|
}\\
&&\le \|v_{\eps \Psi\Phi_i}^{(2)}(\vphi_i(p_\eps))
- v_{\eps \Psi\Phi_i}^{(2)}(\vphi_i(p_\eps'))\|\, |\bfs{\varphi}_i(e_\eps)| 
\\
&&\hphantom{mmmi}\!+\|v_{\eps \Psi\Phi_i}^{(2)}(\vphi_i(p_\eps'))\|\,
|\bfs{\vphi}_i(e_\eps) - \bfs{\vphi}_i(e_\eps')|
\eeas
which can be estimated by applying Lemma 
\ref{kriegllem}. 
The claim now follows as in the proof of Theorem 
\ref{compmfwell} since the constants appearing
in this estimate due to the application of Lemma \ref{kriegllem} are independent of
$\eps$: for each $i$ the domain $\vphi_i(V_i\cap \underline{v_\eps}^{-1}(W))$ of 
$v_{\eps \Psi\Phi_{i}}^{(2)}$ ($1\le i \le r$) 
contains the compact set $\vphi_{i}(\overline{V_{i}'})$.

Finally, let $(v_\eps)_\eps \sim_{vb} (v_\eps')_\eps$ (in fact, it suffices to suppose
$(v_\eps)_\eps \sim_{vb0} (v_\eps')_\eps$). Then using the same techniques
as above, an estimate of the type 
$$
|\bfs{\psi}(v_\eps(e_\eps))- \bfs{\psi}(v'_\eps(e_\eps))|
\le \|v_{\eps \Psi\Phi}^{(2)}(\vphi(p_\eps)) - 
{v'}_{\eps \Psi\Phi}^{(2)}(\vphi(p_\eps))\||\bfs{\vphi}(e_\eps)|
$$
yields that $(v_\eps(e_\eps))_\eps \sim_{vb} (v_\eps'(p_\eps))_\eps$.
\ep

\vskip6pt
Furthermore, a point value characterization based on these concepts is given
in the following result. 

\bt 
\begin{itemize}
\item[(i)] 
Let $u\in \Gamma_{\gs}(X,E)$. Then $u=0$ if and only if $u(\tilde p) = 
0_{\tilde p}$ in $E_c^{\sim_{vb}}$ for all $\tilde p\in \tilde X_c$. 
\item[(ii)] 
Let $(v_\eps)_\eps$, $(v_\eps')_\eps\in 
{\esm}^{\mathrm{VB}}[E,F]$ as in 
Proposition \ref{pvwellvb}\, (ii). 
Then $(v_\eps)_\eps \sim_{vb0} (v_\eps')_\eps$ if and only
if $v(\tilde e) = v'(\tilde e)$ 
for all $\tilde e \in E_c^{\sim_{vb}}$ 
\end{itemize}
\et
\pr 
Necessity is immediate from Proposition \ref{pvwellvb} \, (i). Conversely, suppose that
$u\not=0$ in $\Gamma_{\gs}(X,E)$. Then for some representative $(u_\eps)_\eps$
we have
$$
\exists K\comp X \  \exists m_0\in \N_0 \  \forall k\in \N \  \exists p_k\in K
\  \exists \eps_k < \frac{1}{k}: \quad \|u_{\eps_k}(p_{\eps_k})\|_h
> k \eps_k^{m_0}\,.
$$
Setting $p_\eps = p_{\eps_k}$ for $\eps_{k+1} < \eps \le \eps_k$ we obtain an
element $\tilde p = \cl[(p_\eps)_\eps]$ of $\tilde X_c$ with  
$u(\tilde p) \not=0_{\tilde p}$.

(ii) Necessity follows from the proof of Proposition \ref{pvwellvb}\, (ii).
Conversely, let $(v_\eps)_\eps
\not\sim_{vb0} (v_\eps')$. 
Then either $(\underline{u_\eps})_\eps \not\sim_0 (\underline{v_\eps})_\eps$ holds or 
((ii) for $k=0$) from Definition \ref{homgequ} is violated.  In the first case, 
by Proposition \ref{pvthmfval} there exists $\tilde p \in \tilde X_c$ with 
$\underline{v}(\tilde p) \not= 
\underline{v'}(\tilde p)$. But then $v(0_{\tilde p}) \not= v'(0_{\tilde p})$.
Finally, suppose $\neg (\mbox{(ii) for}\ k=0)$. Then
there exist $m\in\N$, vector bundle charts
$(V,\Phi)$, $(W,\Psi)$ and 
$L\comp V$, $L'\comp W$ allowing the construction of
sequences $\eps_j(<1/j)$, $p_j\in
L\cap \underline{v_{\eps_j}}^{-1}(L')\cap
\underline{{v'}_{\eps_j}}^{-1}(L')$ (for $j=1,2,\dots$) satisfying
\[
\|(v_{\eps_j \Psi \Phi}^{(2)}
- {v'}_{\eps_j \Psi \Phi}^{(2)})(\vphi(p_j))\| > j \eps_j^{m}.
\]
Hence there exists a bounded sequence of vectors $v_j \in \R^{n'}$ such that
$$
\|v_{\eps_j \Psi \Phi}^{(2)}(x_j)v_j
- {v'}_{\eps_j \Psi \Phi}^{(2)}(x_j)v_j\| > j \eps_j^{m}
$$
(with $x_j:=\vphi(p_j)$). Now set
$e_\eps := \Phi^{-1}(x_j,v_j)$ for $\eps_{j+1} < \eps \le \eps_j$
to obtain an element $\tilde e$ of $E_c^{\sim_{vb}}$ with 
$v(\tilde e) \not= v'(\tilde e)$.
\ep

\vskip6pt
Another important transfer process made possible by the above point value concept concerns
the alternative description of smooth tensor fields both as sections and multilinear
mappings and the resulting character of tensor fields as ``pointwise'' objects.
The following result secures this property also in the context of generalized
tensor fields:
\bp \label{gentensorpunktal} 
Let $s \in (\gs)^r_s(X)$ and $\tilde p \in \tilde X_c$. Let $\om_1,\dots,\om_r$
and $\om_1',\dots,\om_r'$ be generalized one-forms with $\om_i(\tilde p) = 
\om_i'(\tilde p)$ ($1\le i \le r$). Let $\xi_1,\dots,\xi_s$
and $\xi_1',\dots,\xi_s'$ be generalized vector fields with $\xi_i(\tilde p) = 
\xi_i'(\tilde p)$ ($1\le i \le s$). Then
$$
(s(\om_1,\dots,\om_r,\xi_1,\dots,\xi_s))(\tilde p)
= (s(\om'_1,\dots,\om'_r,\xi'_1,\dots,\xi'_s))(\tilde p).
$$  
\et
\pr As in the proof of the analogous result for smooth tensor fields it suffices
to show that for $\om \in ({\gs})^0_1(X)$ and 
$\xi\in (\gs)^1_0(X)$ with $\xi(\tilde p)=
0$ we have $(\om(\xi))(\tilde p) = 0$. 
To begin with, let us suppose that there exists a chart $(V,\vphi)$ and $K\comp V$
such that $p_\eps \in K$ for $\eps$ small.

In the local coordinates $(V,\vphi)$ we may write $\xi = 
\sum \xi_j \pa_j$ with $\xi_i \in \gs(V)$. Choose a smooth bump 
function $f$ supported in $V$ such that $f|_K \equiv 1$
in a neighborhood of $K$. Then 
$F:= \cl[(f)_\eps] \in \gs(V)$ and 
$$
F^2\om(\xi) = \om(F^2\xi) = \om(\sum F\xi_j F\pa_j) = \sum F\xi_j \om(F\pa_j)
$$
from which the result follows upon inserting $\tilde p$.

In the general case, let us suppose to the contrary that $\om(\xi)(\tilde p) \not= 0$.
Then there exist
$m_0\in \N$, $\eps_k\to 0$ and representatives $(\om_\eps)_\eps$, $(\xi_\eps)_\eps$, 
$(p_\eps)_\eps$ such that 
\begin{equation}
  \label{tenslocpoint}
|\om_{\eps_k}(\xi_{\eps_k})(p_{\eps_k})|>\eps_k^{m_0} 
\end{equation}
for all $k\in \N$. Without 
loss of generality we may suppose that $p_{\eps_k}$ converges to some $p\in X$.
Choose some chart $(V,\vphi)$ containing $p$ and some $K\comp V$ such that 
$p_{\eps_k}\in K$ for large $k$. Now set $p_\eps' = p_{\eps_k}$ 
and $\xi_\eps' = \xi_{\eps_k}$ for $\eps_{k+1}\le \eps\le \eps_k$. Then we
obtain representatives of elements $\xi'$ of $(\gs)^1_0(X)$ and $\tilde p'$ of
$\tilde X_c$ with $p_\eps' \in K$ for $\eps$ small
such that $\xi'(\tilde p')=0$ but $\om(\tilde \xi)(\tilde p)\not=0$ by
(\ref{tenslocpoint}), a contradiction to what we have already proved in the special
case above.
\ep\ms
{\bf Acknowledgments.} I am greatly indebted to Michael Grosser for several suggestions
that importantly contributed to the final form of the paper. Also, I would like
to thank Roland Steinbauer and Andreas Kriegl for helpful discussions.

\end{document}